\documentstyle[12pt]{article}
   \title{{\bf Application of a ``Jacobi identity'' for vertex
operator algebras to zeta values and differential operators}}
    \author{J. Lepowsky}
    \date{\it Dedicated to the memory of Mosh\'e Flato}
    \begin{document}

    \bibliographystyle{alpha} 
    \maketitle

\begin{abstract}  
We explain how to use a certain new ``Jacobi identity'' for vertex
operator algebras, announced in a previous paper (math.QA/9909178), to
interpret and generalize recent work of S.~Bloch's relating values of
the Riemann zeta function at negative integers with a certain Lie
algebra of operators.
\end{abstract}

    \input amssym.def
    \input amssym  
 
    \newtheorem{rema}{Remark}[section]
    \newtheorem{propo}[rema]{Proposition}
    \newtheorem{theo}[rema]{Theorem}
    \newtheorem{defi}[rema]{Definition}
    \newtheorem{lemma}[rema]{Lemma}
    \newtheorem{corol}[rema]{Corollary}
    \newtheorem{exam}[rema]{Example}

\newcommand{\nordplus}{\mbox{\scriptsize ${+ \atop +}$}}
\newcommand{\nordbullet}{\mbox{\tiny ${\bullet\atop\bullet}$}}

\setcounter{section}{0}
\renewcommand{\theequation}{\thesection.\arabic{equation}}
\renewcommand{\therema}{\thesection.\arabic{rema}}
\setcounter{equation}{0}
\setcounter{rema}{0}

\section{Introduction}

In \cite{L2} we have announced some results relating vertex operator
algebra theory to certain phenomena associated with values of the
Riemann zeta function $\zeta(s)$ at negative integers, and we also
presented expository background.  In the present paper, a
continuation, we announce and sketch the proofs of more of these
results.  The goal is to use vertex operator algebra theory to
interpret the classical formal equality
\begin{equation}\label{zeta-s}
\sum_{n>0}n^s = \zeta(-s), \;\; s=0,1,2,\dots,
\end{equation}
where the left-hand side is the divergent sum and the right-hand side
is the zeta function (analytically continued), and to exploit the
interpretation in the setting of vertex operator algebra theory.  More
specifically, we wanted to ``explain'' some recent work \cite{Bl} of
S. Bloch's involving zeta values and central extensions of Lie
algebras of differential operators.  Below we shall summarize the
contents of \cite{L2} and make some additional comments, and then in
Section 3 we shall present the new results.  These deal with the
passage {}from the general ``Jacobi identity'' announced in \cite{L2},
Theorem 4.2 (recalled in Theorem \ref{newjacobithm} below), to the
more special formula, Theorem 3.1 of \cite{L2} (recalled in Theorem
\ref{theoremforLbar} below).  This in turn recovered the main formulas
in \cite{Bl}, expressing the structure of a Lie algebra, constructed
using zeta-function values {}from certain vertex operators, providing
a central extension of a Lie algebra of differential operators.  The
work \cite{L3} contains details and related results.

While it might help the reader of the present paper to have the paper
\cite{L2} available to refer to, we shall review some background and
motivation in this extended Introduction and in Section 2 below, to
make this paper more self-contained.  The famous classical ``formula''
\begin{equation}\label{1+}
1 + 2 + 3 + \cdots = - \frac{1}{12},
\end{equation}
which has the rigorous meaning
\begin{equation}\label{zeta-1}
\zeta(-1) = - \frac{1}{12},
\end{equation}
is intimately related to the regularizing of certain infinities in
conformal field theory.  In \cite{L2} we announced some general
principles of vertex operator algebra theory that elucidate the
passage {}from the unrigorous but suggestive formula (\ref{1+}) to
formula (\ref{zeta-1}), and the generalization (\ref{zeta-s}), in the
process ``explaining'' some results in \cite{Bl}.  Foundational
notions of vertex operator algebra theory, and ``formal calculus,'' as
presented in \cite{FLM} and \cite{FHL}, enter in an essential way.

Our main themes involve: the ubiquity of generating functions---the
introduction of new formal variables and generating functions in order
to render complicated things easier, more natural, and at the same
time, much more general, as in the corresponding parts of \cite{FLM};
the use of commuting formal variables rather than complex variables
because they provide the most natural way to handle the
doubly-infinite series such as $\delta (x) = \sum_{n \in {\Bbb Z}}x^n$
that pervade the natural formulations and proofs in the subject; the
exploitation of the formal exponential of the differential operator
$x{\frac{d}{dx}}$ as a formal change-of-variables automorphism (again
as in \cite{FLM}); the formulation of Euler's interpretation of the
divergent series (\ref{zeta-s}) by means of the operator product
expansion in conformal field theory; and the consideration of central
extensions of Lie algebras of differential operators in terms of the
very general context of what we termed the ``Jacobi identity''
\cite{FLM} for vertex operator algebras.

Now we sketch the relevant classical background on the Virasoro
algebra and related matters.  Consider the Lie algebra
\begin{equation}
{\frak {d}} = \mbox{\rm Der}\ \Bbb C [t,t^{-1}]
\end{equation}
of formal vector fields on the circle, with basis
$\{t^n D | n \in \Bbb Z\}$, where
\begin{equation}\label{D}
D = D_t = t \frac{d}{dt},
\end{equation}
and recall the Virasoro algebra ${\frak {v}}$, the central extension
\begin{equation}\label{Vir}
0 \rightarrow \Bbb C c \rightarrow {\frak {v}} \rightarrow {\frak {d}}
\rightarrow 0,
\end{equation}
where ${\frak {v}}$ has basis $\{L(n) | n \in \Bbb Z \}$ together with
a central element $c$; the bracket relations are given by
\begin{equation}\label{Virbrackets}
[L(m),L(n)] = (m-n)L(m+n) + \frac{1}{12} (m^3 - m) \delta_{m+n,0} c
\end{equation}
and $L(n)$ maps to $-t^n D$ in (\ref{Vir}).  This Lie algebra is
naturally $\Bbb Z$-graded, with ${\rm deg}\ L(n) = n$ and ${\rm deg}\
c = 0$.  It has the following well-known realization: Start with the
Heisenberg Lie algebra with basis consisting of the symbols $h(n)$ for
$n \in \Bbb Z$, $n \neq 0$ and a central element $1$, with the bracket
relations
\begin{equation}\label{Heisbrackets}
[h(m),h(n)] = m \delta_{m+n,0} 1.
\end{equation}
For convenience we adjoin an additional central basis element $h(0)$,
so that the relations (\ref{Heisbrackets}) hold for all $m,n \in
\Bbb Z$.  This Lie algebra acts irreducibly on the polynomial algebra
\begin{equation}\label{S}
S = \Bbb C [h(-1),h(-2),h(-3),\dots]
\end{equation}
as follows: For $n<0$, $h(n)$ acts by multiplication; for $n>0$,
$h(n)$ acts as $n \frac{\partial}{\partial h(-n)}$; $h(0)$ acts as 0;
and $1$ acts as the identity operator.  Then ${\frak {v}}$ acts on $S$
via:
\begin{equation}\label{c}
c \mapsto 1,
\end{equation}
\begin{equation}\label{Ln}
L(n) \mapsto \frac{1}{2} \sum_{j \in \Bbb Z} h(j)h(n-j) \ \ \mbox{for}
\ n \neq 0,
\end{equation}
\begin{equation}\label{L0}
L(0) \mapsto \frac{1}{2} \sum_{j \in \Bbb Z} h(-|j|)h(|j|).
\end{equation}

In the case of $L(0)$, the absolute values make the operator well
defined, while for $n \neq 0$ the operator is well defined as it
stands, since $[h(j),h(n-j)] = 0$; thus the operators (\ref{Ln}) and
(\ref{L0}) are in ``normal-ordered form,'' that is, the $h(n)$ for
$n>0$ act to the right of the $h(n)$ for $n<0$.  Using colons to
denote normal ordering, we thus have
\begin{equation}\label{Lnnormal}
L(n) \mapsto \frac{1}{2} \nordbullet\sum_{j \in \Bbb Z}
h(j)h(n-j)\nordbullet
\end{equation}
for all $n \in \Bbb Z$.  It is well known that the operators
(\ref{Lnnormal}) indeed satisfy the bracket relations
(\ref{Virbrackets}).  (This exercise and the related constructions are
presented in \cite{FLM}, for example, where the standard
generalization of this construction of ${\frak {v}}$ using a
Heisenberg algebra based on a finite-dimensional space of operators
$h(n)$ for each $n$ is also carried out.)  In vertex operator algebra
theory and conformal field theory it is standard procedure to embed
operators such as $h(n)$ and $L(n)$ into generating functions and to
compute with these generating functions, using a formal calculus.

The space $S$ is naturally $\Bbb Z$-graded, with ${\rm deg}\ h(j) = j$
for $j < 0$, and $S$ is graded as a ${\frak {v}}$-module.  It is
appropriate to use the negative of this grading, that is, to define a
new grading (by ``conformal weights'') on the space $S$ by the rule
${\rm wt}\ h(-j) = j$ for $j > 0$; for each $n \geq 0$, the
homogeneous subspace of $S$ of weight $n$ coincides with the
eigenspace of the operator $L(0)$ with eigenvalue $n$.  For $n \in
{\Bbb Z}$ (or $n \geq 0$) let $S_n$ be the homogeneous subspace of $S$
of weight $n$, and consider the formal power series
\begin{equation}
\mbox{dim}_* S = \sum_{n \geq 0} (\mbox{dim}\;S_n) q^n
\end{equation}
(the ``graded dimension'' of the graded space $S$).  Clearly,
\begin{equation}
\mbox{dim}_* S = \prod_{n > 0} (1 - q^n)^{-1}.
\end{equation}

As is well known, removing the normal ordering in the definition of
$L(0)$ introduces an infinity which formally equals $\frac{1}{2}
\zeta(-1)$, since the unrigorous expression
\begin{equation}\label{Lbar(0)}
{\bar L}(0) = \frac{1}{2} \sum_{j \in \Bbb Z} h(-j)h(j)
\end{equation}
formally equals (by (\ref{Heisbrackets}))
\begin{equation}\label{Lbar(0)2}
L(0) + \frac{1}{2} (1 + 2 + 3 + \cdots),
\end{equation}
which itself formally equals
\begin{equation}\label{Lbar(0)3}
L(0) + \frac{1}{2} \zeta(-1) = L(0) - \frac{1}{24}.
\end{equation}
Rigorizing ${\bar L}(0)$ by defining it as
\begin{equation}\label{Lbar(0)rig}
{\bar L}(0) = L(0) + \frac{1}{2} \zeta(-1),
\end{equation}
we set
\begin{equation}\label{Lbar(n)}
{\bar L}(n) = L(n) \ \ \mbox{for} \ n \neq 0,
\end{equation}
to get a new basis of ${\frak {v}}$.  (We are identifying the elements
of ${\frak {v}}$ with operators on $S$.)  The brackets become:
\begin{equation}\label{newVirbrackets}
[{\bar L}(m),{\bar L}(n)] = (m-n){\bar L}(m+n) + \frac{1}{12} m^3
\delta_{m+n,0};
\end{equation}
that is, $m^3 - m$ in (\ref{Virbrackets}) has become the pure monomial
$m^3$.

It is a fundamental and well-known fact that this formal removal of
the normal ordering gives rise to modular transformation properties:
We define a new grading on $S$ by using the eigenvalues of ${\bar
L}(0)$ in place of $L(0)$, so that the grading of $S$ is ``shifted''
{}from the previous grading by conformal weights by the subtraction of
$\frac{1}{24}$ {}from the weights.  Letting $\chi (S)$ be the
corresponding graded dimension, we have
\begin{equation}\label{chiS}
\chi (S) = \frac{1}{\eta (q)},
\end{equation}
where
\begin{equation}
\eta (q) = q^{\frac{1}{24}} \prod_{n > 0} (1 - q^n). 
\end{equation}
With the formal variable $q$ replaced by $e^{2 \pi i \tau}$, $\tau$ in
the upper half-plane, $\eta (q)$, Dedekind's eta-function, has
important (classical) modular transformation properties, unlike
$\prod_{n > 0} (1 - q^n)$.

In \cite{Bl}, Bloch considered the larger Lie algebra of formal
differential operators, spanned by
\begin{equation}
\{ t^n D^m | n \in \Bbb Z, \ m \geq 0 \}
\end{equation}
or more precisely, we restrict to $m > 0$ and further, to the Lie
subalgebra ${\cal D}^+$, containing $\frak {d}$, spanned by the
differential operators of the form $D^r (t^nD) D^r$ for $r \geq 0, \ n
\in \Bbb Z$.  Then we can construct a central extension of ${\cal
D}^+$ using generalizations of the operators (\ref{Lnnormal}):
\begin{equation}\label{Lrn}
L^{(r)}(n) = \frac{1}{2} \sum_{j \in \Bbb Z} j^r h(j) (n-j)^r h(n-j) \
\ \mbox{for} \ n \neq 0,
\end{equation}
\begin{equation}\label{Lr0}
L^{(r)}(0) = \frac{1}{2} \sum_{j \in \Bbb Z} (-j)^r h(-|j|) j^r
h(|j|),
\end{equation}
that is,
\begin{equation}\label{Lrnnormal}
L^{(r)}(n) = \nordbullet\frac{1}{2} \sum_{j \in \Bbb Z} j^r h(j) (n-j)^r
h(n-j)\nordbullet
\end{equation}
for $n \in \Bbb Z$.  These operators provide \cite{Bl} a central
extension of ${\cal D}^+$ such that
\begin{equation}
L^{(r)}(n) \mapsto (-1)^{r+1} D^r (t^nD) D^r.
\end{equation}

A central point of \cite{Bl} is that the formal removal of the
normal-ordering procedure in the definition (\ref{Lr0}) of
$L^{(r)}(0)$ adds the infinity $(-1)^r \frac{1}{2} \zeta(-2r-1) =$ ``
$\sum_{n>0} n^{2r+1}$'' (generalizing
(\ref{Lbar(0)})--(\ref{Lbar(0)3})), and if we correspondingly define
\begin{equation}\label{lbarr(0)}
{\bar L}^{(r)} (0) = L^{(r)} (0) + (-1)^r \frac{1}{2} \zeta(-2r-1)
\end{equation}
and ${\bar L}^{(r)} (n) = L^{(r)} (n)$ for $n \neq 0$ (generalizing
(\ref{Lbar(0)rig}) and (\ref{Lbar(n)})), the commutators simplify in a
remarkable way: The complicated polynomial in the scalar term of
$[{\bar L}^{(r)} (m),{\bar L}^{(s)} (-m)]$ reduces to a pure monomial
in $m$, by analogy with, and generalizing, the passage {}from $m^3 -
m$ to $m^3$ in (\ref{newVirbrackets}); see \cite{Bl} for the formulas
and further results.  (See also \cite{PRS} and \cite{FS} for
treatments of essentially the same phenomena in physical contexts.)

In \cite{L2} two layers of ``explanation'' and generalization of the
results of \cite{Bl} were presented.  It was recalled that for $k >
1$,
\begin{equation}\label{zetaB}
\zeta(-k+1) = - \frac{B_{k}}{k},
\end{equation}
where the $B_{k}$ are the Bernoulli numbers, defined by the generating
function
\begin{equation}\label{Ber}
\frac{x}{e^x - 1} = \sum_{k \geq 0} \frac{B_{k}}{k!} x^k,
\end{equation}
where $x$ is a formal variable, and (a variant of) Euler's heurstic
interpretation of (\ref{zetaB}) in terms of the divergent series
(\ref{zeta-s}) was also recalled.  The ``first layer of explanation''
proceeded as follows:

Using a formal variable $x$, we form the generating functions
\begin{equation}\label{h(x)}
h(x) = \sum_{n \in \Bbb Z} h(n)x^{-n}
\end{equation}
and
\begin{equation}
L^{(r)}(x) = \sum_{n \in \Bbb Z} L^{(r)}(n)x^{-n},
\end{equation}
and using $D_x$ to denote the operator $x \frac{d}{dx}$ (as in
(\ref{D})), we observe that
\begin{equation}\label{Lrx}
L^{(r)}(x) = {\frac{1}{2}}\nordbullet(D_x^r h(x))^2\nordbullet,
\end{equation}
where the colons, as always, denote normal ordering.  (For other
purposes, other versions of these generating functions are used, in
particular, $h(x) = \sum_{n \in \Bbb Z} h(n)x^{-n-1}$ in place of
(\ref{h(x)}); see (\ref{Yh(-1)}) below.)

Next we introduce suitable generating functions over the number of
{\it derivatives}, and we use the formal multiplicative analogue
\begin{equation}\label{infinitdil}
e^{y D_x} f(x) = f(e^y x)
\end{equation}
of the formal Taylor theorem
\begin{equation}\label{Taylor}
e^{y \frac{d}{dx}} f(x) = f(x+y),
\end{equation}
where $f(x)$ is an arbitrary formal series of the form $\sum_n a_n
x^n$, $n$ is allowed to range over something very general, like $\Bbb
Z$ or even $\Bbb C$, say, and the $a_n$ lie in a fixed vector space
(cf. \cite{FLM}, Proposition 8.3.1).  Although $\nordbullet(D_x^r
h(x))^2\nordbullet$ (recall (\ref{Lrx})) is hard to put into a
``good'' generating function over $r$, we make the problem easier by
making it more general: Consider independently many derivatives on
each of the two factors $h(x)$ in $\nordbullet h(x)^2\nordbullet$, use
two new independent formal variables $y_1$ and $y_2$, and form the
generating function
\begin{equation}\label{Ly1y2}
L^{(y_1,y_2)}(x) = {\frac{1}{2}}
\nordbullet(e^{y_1 D_x} h(x))(e^{y_2 D_x}
h(x))\nordbullet = {\frac{1}{2}}
\nordbullet h(e^{y_1} x)h(e^{y_2} x)\nordbullet
\end{equation}
(where we use (\ref{infinitdil})), so that $L^{(r)}(x)$ is a
``diagonal piece'' of this generating function.  Using formal vertex
operator calculus techniques, we can calculate
\begin{equation}\label{bracketofquadratics}
[\nordbullet h(e^{y_1} x_1)h(e^{y_2} x_1)\nordbullet,\nordbullet
h(e^{y_3} x_2)h(e^{y_4} x_2)\nordbullet].
\end{equation}

Now the formal expression $h(e^{y_1} x)h(e^{y_2} x)$ is not rigorous,
as we see by (for example) trying to compute the constant term in the
variables $y_1$ and $y_2$ in this expression; the failure of this
expression to be defined in fact corresponds exactly to the occurrence
of formal sums like $\sum_{n>0} n^r$ with $r>0$, as we have been
discussing.  However, we have
\begin{equation}\label{hx1hx2}
h(x_1)h(x_2) = \nordbullet h(x_1)h(x_2)\nordbullet + x_2
\frac{\partial}{\partial x_2}
\frac{1}{1 - x_2 / x_1}
\end{equation}
and it follows that
\begin{equation}\label{hex1hex2}
h(e^{y_1} x_1)h(e^{y_2} x_2) = \nordbullet h(e^{y_1} x_1)h(e^{y_2}
x_2)\nordbullet + x_2
\frac{\partial}{\partial x_2} \frac{1}{1 - e^{y_2} x_2 / e^{y_1}
x_1};
\end{equation}
note that $x_2\frac{\partial}{\partial x_2}$ can be replaced by $-
\frac{\partial}{\partial y_1}$ in the last expression.  The expression
$\frac{1}{1 - e^{y_2} x_2 / e^{y_1} x_1}$ came {}from, and is, a
geometric series expansion (recall (\ref{hx1hx2})).

If we try to set $x_1 = x_2 \ (= x)$ in (\ref{hex1hex2}), the result
is unrigorous on the left-hand side, as we have pointed out, {\it but
the result has rigorous meaning on the right-hand side}, because the
normal-ordered product $\nordbullet h(e^{y_1} x)h(e^{y_2}
x)\nordbullet$ is certainly well defined, and the expression $-
\frac{\partial}{\partial y_1} \frac{1}{1 - e^{-y_1 + y_2}}$ can be
interpreted rigorously as in (\ref{Ber}); more precisely, we take
$\frac{1}{1 - e^{-y_1 + y_2}}$ to mean the formal (Laurent) series in
$y_1$ and $y_2$ of the shape
\begin{equation}\label{defof1/1-e}
\frac{1}{1 - e^{-y_1 + y_2}} = (y_1 - y_2)^{-1}F(y_1,y_2),
\end{equation}
where $(y_1 - y_2)^{-1}$ is understood as the binomial expansion
(geometric series) in nonnegative powers of $y_2$ and $F(y_1,y_2)$ is
an (obvious) formal power series in (nonnegative powers of) $y_1$ and
$y_2$.  This motivates us to define a new ``normal-ordering''
procedure
\begin{equation}\label{hexhex}
\nordplus h(e^{y_1} x)h(e^{y_2} x)\nordplus = 
\nordbullet h(e^{y_1} x)h(e^{y_2} x)\nordbullet -
\frac{\partial}{\partial y_1} \frac{1}{1 - e^{-y_1 + y_2}},
\end{equation}
with the last part of the right-hand side being understood as we just
indicated.  {\it This ``rigorization'' of the undefined formal
expression $h(e^{y_1} x)h(e^{y_2} x)$ corresponds exactly to Euler's
heuristic interpretation of (\ref{zetaB}) discussed in \cite{L2}.}
This gives us a natural ``explanation'' of the zeta-function-modified
operators defined in (\ref{lbarr(0)}): We use (\ref{hexhex}) to define
the following analogues of the operators (\ref{Ly1y2}):
\begin{equation}\label{Lbary1y2}
{\bar L}^{(y_1,y_2)}(x) = 
{\frac{1}{2}}
\nordplus h(e^{y_1} x)h(e^{y_2} x)\nordplus,
\end{equation}
and the operator ${\bar L}^{(r)} (n)$ is exactly $(r!)^2$ times the
coefficient of $y_1^r y_2^r x_0^{-n}$ in (\ref{Lbary1y2}); the
significant case is the case $n = 0$.

With the new normal ordering (\ref{hexhex}) replacing the old one,
remarkable cancellation occurs in the commutator
(\ref{bracketofquadratics}), and here is the main result of the
``first layer of explanation'' of Bloch's formula for $[{\bar L}^{(r)}
(m),{\bar L}^{(s)} (n)]$ in \cite{Bl}, in a somewhat generalized form:
\begin{theo} \cite{L2} \label{theoremforLbar}
With the formal delta-function Laurent series $\delta(x)$ defined as
\begin{equation}\label{delta}
\delta(x) = \sum_{n \in {\Bbb Z}} x^n,
\end{equation}
and with independent commuting formal variables as indicated, we have:
\begin{eqnarray}\label{Lbarbrackets}
\lefteqn{[{\bar L}^{(y_1,y_2)}(x_1),{\bar L}^{(y_3,y_4)}(x_2)]}\nonumber\\
&&= - {\frac{1}{2}} \frac{\partial}{\partial y_1}	
\biggl({\bar L}^{(-y_1+y_2+y_3,y_4)}(x_2)
\delta \left({\frac{e^{y_1}x_1}{e^{y_3}x_2}}\right)\nonumber\\
&&\hspace{2em} + {\bar L}^{(-y_1+y_2+y_4,y_3)}(x_2)
\delta \left({\frac{e^{y_1}x_1}{e^{y_4}x_2}}\right)\biggr)\nonumber\\
&&\hspace{2em} - {\frac{1}{2}} \frac{\partial}{\partial y_2}	
\biggl({\bar L}^{(y_1-y_2+y_3,y_4)}(x_2)
\delta \left({\frac{e^{y_2}x_1}{e^{y_3}x_2}}\right)\nonumber\\
&&\hspace{2em} + {\bar L}^{(y_1-y_2+y_4,y_3)}(x_2)
\delta \left({\frac{e^{y_2}x_1}{e^{y_4}x_2}}\right)\biggr).
\end{eqnarray}
\end{theo}

As discussed in \cite{L2}, the pure monomials in $m$ that we set out
to explain now emerge completely naturally, since for example the
delta-function expression $\delta \left(e^{y_1}x_1 / e^{y_3}x_2
\right)$ can be written as $e^{y_1 D_{x_1}}e^{y_3 D_{x_2}}\delta
\left(x_1 / x_2 \right)$, and when we extract and equate the
coefficients of the monomials in the variables $y_1^r y_2^r y_3^s
y_4^s$ on the two sides of (\ref{Lbarbrackets}), we get expressions
like $(D^j \delta)\left(x_1 / x_2 \right)$, whose expansion, in turn,
in powers of $x_1$ and $x_2$ clearly yields a pure monomial analogous
to and generalizing the expression $m^3$ in (\ref{newVirbrackets}).
(All of these considerations hold equally well in the more general
situation where we start with a Heisenberg algebra based on a
finite-dimensional space rather than a one-dimensional space.)

In the next section we recall {}from \cite{L2} how these
considerations can actually be understood in the much greater
generality of vertex operator algebra theory in general (the ``second
explanation and generalization'').

I am very grateful to Spencer Bloch for informing me about his work
and for many valuable discussions.

This work was partially supported by NSF grants DMS-9401851 and
DMS-9701150.

\renewcommand{\theequation}{\thesection.\arabic{equation}}
\renewcommand{\therema}{\thesection.\arabic{rema}}
\setcounter{equation}{0}
\setcounter{rema}{0}

\section{The general principles}

At this point we recall the definition of the notion of vertex
operator algebra {}from \cite{FLM}, \cite{FHL}.  This is a variant of
Borcherds' notion of vertex algebra \cite{Bo} and is based on the
``Jacobi identity'' as formulated in \cite{FLM} and \cite{FHL}; we
need this formal-variable formulation in order to express our results
in a natural way.  We continue to use commuting formal variables $x,
x_0, x_1, x_2$ and several other variables.  Recall the formal
delta-function Laurent series $\delta(x)$ (\ref{delta}).
	
\begin{defi}\label{VOA}
{\rm A {\it vertex operator algebra} $(V, Y, {\bf 1}, \omega)$, or simply $V$
(over ${\Bbb C}$), is a ${\Bbb
Z}$-graded vector space (graded by {\it weights})
\begin{equation}
V=\coprod_{n\in {\Bbb Z}}V_{(n)}; \ \mbox{\rm for}\ v\in V_{(n)},\;n=\mbox{\rm wt}\ v;
\end{equation}
such that
\begin{equation}
\mbox{\rm dim }V_{(n)}<\infty\;\;\mbox{\rm for}\; n \in {\Bbb Z},
\end{equation}
\begin{equation}
V_{(n)}=0\;\;\mbox{\rm for} \;n\; \mbox{\rm sufficiently small},
\end{equation}
equipped with a linear map  $V\otimes V\to V[[x, x^{-1}]]$, or
equivalently,
\begin{eqnarray}
V&\to&(\mbox{\rm End}\; V)[[x, x^{-1}]]\nonumber \\
v&\mapsto& Y(v, x)={\displaystyle \sum_{n\in{\Bbb Z}}}v_{n}x^{-n-1}
\;\;(\mbox{\rm where}\; v_{n}\in
\mbox{\rm End} \;V),
\end{eqnarray}
$Y(v, x)$ denoting the {\it vertex operator associated with} $v$, and
equipped also with two distinguished homogeneous vectors ${\bf 1}\in
V_{(0)}$ (the {\it vacuum}) and $\omega \in V_{(2)}$. The following
conditions are assumed for $u, v \in V$: the {\it lower truncation
condition} holds:
\begin{equation}
u_{n}v=0\;\;\mbox{\rm for}\;n\; \mbox{\rm sufficiently large}
\end{equation}
(or equivalently, $Y(u,x)v$ involves only finitely many negative
powers of $x$);
\begin{equation}
Y({\bf 1}, x)=1\;\; (1\;\mbox{\rm on the right being the identity
operator});
\end{equation}
the {\it creation property} holds:
\begin{equation}
Y(v, x){\bf 1} \in V[[x]]\;\;\mbox{\rm and}\;\;\lim_{x\rightarrow
0}Y(v, x){\bf 1}=v
\end{equation}
(that is, $Y(v, x){\bf 1}$ involves only nonnegative integral powers
of $x$ and the constant term is $v$); with binomial expressions
understood to be expanded in nonnegative powers of the second
variable, the {\it Jacobi identity} (the main axiom) holds:
\begin{eqnarray}\label{jacobi}
&x_{0}^{-1}\delta
\left({\displaystyle\frac{x_{1}-x_{2}}{x_{0}}}\right)Y(u, x_{1})Y(v,
x_{2})-x_{0}^{-1} \delta
\left({\displaystyle\frac{x_{2}-x_{1}}{-x_{0}}}\right)Y(v, x_{2})Y(u,
x_{1})&\nonumber \\ &=x_{2}^{-1} \delta
\left({\displaystyle\frac{x_{1}-x_{0}}{x_{2}}}\right)Y(Y(u, x_{0})v,
x_{2})&
\end{eqnarray}
(note that when each expression in (\ref{jacobi}) is applied to any
element of $V$, the coefficient of each monomial in the formal
variables is a finite sum; on the right-hand side, the notation
$Y(\cdot, x_{2})$ is understood to be extended in the obvious way to
$V[[x_{0}, x^{-1}_{0}]]$); the Virasoro algebra relations hold (acting
on $V$):
\begin{equation}
[L(m), L(n)]=(m-n)L(m+n)+{\displaystyle\frac{1}{12}}
(m^{3}-m)\delta_{n+m,0}({\rm rank}\;V)1
\end{equation}
for $m, n \in {\Bbb Z}$, where
\begin{equation}
L(n)=\omega _{n+1}\;\; \mbox{\rm for} \;n\in{\Bbb Z}, \;\;{\rm
i.e.},\;\;Y(\omega, x)=\sum_{n\in{\Bbb Z}}L(n)x^{-n-2}
\end{equation}
and 
\begin{eqnarray}
&{\rm rank}\;V\in {\Bbb C};&\\ 
&L(0)v=nv=(\mbox{\rm wt}\ v)v\;\;\mbox{\rm for}\;n \in {\Bbb
Z}\;\mbox{\rm and}\;v\in V_{(n)};&\\ 
&{\displaystyle \frac{d}{dx}}Y(v,x)=Y(L(-1)v, x)&
\end{eqnarray}
(the {\it  $L(-1)$-derivative property}).}
\end{defi}

In the presence of simpler axioms, the Jacobi identity in the
definition is equivalent to a suitably-formulated ``commutativity''
relation and a suitably formulated ``associativity relation'' (recall
\cite{FLM}, \cite{FHL}).  The commutativity condition asserts that for
$u,v \in V$,
\begin{equation}
Y(u,x_1)Y(v,x_2) \sim Y(v,x_2)Y(u,x_1),
\end{equation}
where ``$\sim$'' denotes equality up to a suitable kind of generalized
analytic continuation, and the associativity condition asserts that
\begin{equation}\label{assoc}
Y(u,x_1)Y(v,x_2) \sim Y(Y(u,x_1 - x_2)v,x_2),
\end{equation}
where the right-hand side and the generalized analytic continuation
have to be understood in suitable ways (see \cite{FLM} and \cite{FHL}
and cf. \cite{BPZ} and \cite{G}); the right-hand side of (\ref{assoc})
is {\it not} a well-defined formal series in $x_1$ and $x_2$.

In \cite{L2} we reformulated the associativity relation in the
following way: Formally replacing $x_1$ by $e^{y}x_2$ in
(\ref{assoc}), we find (formally and unrigorously) that
\begin{equation}\label{assoc2}
Y(u,e^{y}x_2)Y(v,x_2) \sim Y(Y(u,(e^{y}-1)x_2)v,x_2).
\end{equation}
While the left-hand side of (\ref{assoc2}) is not a well-defined
formal series in the formal variables $y$ and $x_2$, the right-hand
side of (\ref{assoc2}) {\it is} in fact a well-defined formal series
in these formal variables.  By replacing $x_1$ by $e^{y}x_2$ we have
made the right-hand side of (\ref{assoc2}) rigorous (and the left-hand
side unrigorous).  Next, instead of the vertex operators $Y(v,x)$, we
want the modified vertex operators defined for homogeneous elements $v
\in V$ by:
\begin{equation}
X(v,x) = x^{{\rm wt}\;v}Y(v,x) = Y(x^{L(0)}v,x),
\end{equation}
as in \cite{FLM}, formula (8.5.27) (recall that $L(0)$-eigenvalues
define the grading of $V$); the formula $X(v,x) = Y(x^{L(0)}v,x)$
works for {\it all} $v \in V$ (not necessarily homogeneous).

Using this we formally obtain {}from (\ref{assoc2}) the formal relation
\begin{equation}\label{assoc3}
X(u,e^{y}x_2)X(v,x_2) \sim X(Y[u,y]v,x_2),
\end{equation}
where $Y[u,y]$ is the operator defined in \cite{Z1}, \cite{Z2} as
follows:
\begin{equation}\label{Ybracket}
Y[u,y] = Y(e^{yL(0)}u,e^{y}-1);
\end{equation}
the right-hand side of (\ref{assoc3}) is still well defined (and the
left-hand side still not well defined).  By Zhu's change-of-variables
theorem in \cite{Z1}, \cite{Z2} (see \cite{L1} and \cite{H1},
\cite{H2} for different treatments of this theorem), $x \mapsto
Y[u,x]$ defines a new vertex operator algebra structure on the same
vector space $V$ under suitable conditions; in \cite{L1}, a simple
formal-variable proof of the Jacobi identity, which we shall need
below, for these operators is given.

The formal relation (\ref{assoc3}) generalizes to products of several
operators, as follows:
\begin{eqnarray}\label{several}
&X(v_1,e^{y_1}x)X(v_2,e^{y_2}x) \cdots X(v_n,e^{y_n}x)X(v_{n+1},x)&
\nonumber \\
&\sim X(Y[v_1,y_1]Y[v_2,y_2] \cdots Y[v_n,y_n]v_{n+1},x).&
\end{eqnarray}
We can also multiply the variable $x$ in this relation by an
exponential, and we get such formal relations as:
\begin{equation}\label{modified}
X(u,e^{y_1}x)X(v,e^{y_2}x) \sim X(Y[u,y_1-y_2]v,e^{y_2}x),
\end{equation}
which should be compared with (\ref{assoc}).

As discussed in \cite{L2}, the point is that very special cases of the
formal relation (\ref{assoc3}) precisely ``explain'' the classical
formal relation (\ref{zeta-s}) in the setting of vertex operator
algebra theory.  Specifically, the polynomial algebra $S$ (recall
(\ref{S})) carries a canonical vertex operator algebra structure of
rank 1 with vacuum vector ${\bf 1}$ equal to $1 \in S$, with the
operators $L(n)$ agreeing with the operators defined in (\ref{Ln}),
(\ref{L0}), with
\begin{equation}\label{Yh(-1)}
Y(h(-1),x) = x^{-1}h(x) = \sum_{n \in \Bbb Z} h(n)x^{-n-1}
\end{equation}
(recall (\ref{h(x)})) and with $\omega = \frac{1}{2} (h(-1))^{2} \in
S$.  Setting
\begin{equation}
v_0 =  h(-1) \in S,
\end{equation}
we thus have
\begin{equation}
X(v_0,x) =  X(h(-1),x) = h(x).
\end{equation}
Next, using this we state a precise formula that equates the new
normal-ordering procedure (\ref{hexhex}) (which in turn interpreted
Bloch's zeta-function-modified operators and the classical formal
relation (\ref{zeta-s})) with its current, still more conceptual,
formulation:
\begin{equation}\label{newnormordrelation}
{\nordplus}h(e^{y}x)h(e^{w}x){\nordplus} = X(Y[v_{0},y-w]v_{0},e^{w}x),
\end{equation}
or equivalently,
\begin{equation}
{\nordplus}h(e^{y+w}x)h(e^{w}x){\nordplus} = X(Y[v_{0},y]v_{0},e^{w}x).
\end{equation}
{\it Moreover, the very general formal relation (\ref{modified}),
applied to the particular parameters in (\ref{newnormordrelation}),
amounts precisely to the formal relation
\begin{equation}\label{hhsimnordplushh}
h(e^{y}x)h(e^{w}x) \sim {\nordplus}h(e^{y}x)h(e^{w}x){\nordplus},
\end{equation}
which was in turn the rigorization of the undefined formal expression
$h(e^{y}x)h(e^{w}x)$ by means of the generating function of the
Bernoulli numbers (recall (\ref{hexhex}).}

Thus we have interpreted (\ref{zeta-s}) and (\ref{hexhex}) by means of
a very special case of a very general picture.  We still need to
``explain'' the bracket relation (\ref{Lbarbrackets}) {}from this
point of view.  (As we have mentioned, everything works for a
Heisenberg Lie algebra based, more generally, on a finite-dimensional
space.)

To do this, we now see that we need to compute in a conceptual way the
commutators of certain expressions of the type $X(Y[u,y]v,x)$.  But
the most natural thing to do (just as was the case in the analogous
contexts in \cite{FLM}, for instance) is to seek a stronger
``Jacobi-type identity,'' analogous to (\ref{jacobi}), and to derive
the desired commutators {}from it.  The main result announced in
\cite{L2} (Theorem 4.2) was in fact just such an identity, for
operators of the simpler type $X(v,x)$ rather than of the type
$X(Y[u,y]v,x)$:

\begin{theo}\label{newjacobithm}
In any vertex operator algebra $V$, for $u,v \in V$ we have:
\begin{eqnarray}\label{newjacobiiden}
&x_{0}^{-1}\delta
\left(e^{y_{21}}{\displaystyle\frac{x_{1}}{x_{0}}}\right)X(u, x_{1})X(v,
x_{2})-x_{0}^{-1} \delta
\left(-e^{y_{12}}{\displaystyle\frac{x_{2}}{x_{0}}}\right)X(v, x_{2})X(u,
x_{1})&\nonumber \\ &=x_{2}^{-1} \delta
\left(e^{-y_{01}}{\displaystyle\frac{x_{1}}{x_{2}}}\right)X(Y[u, y_{01}]v,
x_{2}),&
\end{eqnarray}
where
\begin{equation}\label{yij}
y_{21} = \log \left(1-{\displaystyle\frac{x_{2}}{x_{1}}}\right),\;\;
y_{12} = \log \left(1-{\displaystyle\frac{x_{1}}{x_{2}}}\right),\;\;
y_{01} = -\log \left(1-{\displaystyle\frac{x_{0}}{x_{1}}}\right),\nonumber
\end{equation}
$\log$ denoting the logarithmic formal series.
\end{theo}

The proof is not difficult, starting {}from the Jacobi identity
(\ref{jacobi}).  What was most interesting was that an identity such
as this exists at all, with its strong parallels with the Jacobi
identity (\ref{jacobi}) itself.  The three formal variables defined by
the logarithmic formal series (\ref{yij}) can to a certain extent be
viewed as independent formal variables in their own right; along with
symmetry considerations, this is why we choose to write the
delta-function expressions in (\ref{jacobi}) in this way.

Formula (\ref{newjacobiiden}) can be viewed as solving a problem
implicit in formula (8.8.43) (Corollary 8.8.19) of \cite{FLM}: The
left-hand side of that formula involved the expansion coefficients
with respect to $x_0$, $x_1$ and $x_2$ of the left-hand side of
(\ref{newjacobiiden}), but we had been unable to put the right-hand
side of formula (8.8.43) of \cite{FLM} into an elegant
generating-function form.  Theorem \ref{newjacobithm} solves this
problem, and in the process, in fact relates the left-hand side to
zeta-function values, as explained in the present work.

In \cite{L2} we emphasized the philosophy to ``always use generating
functions,'' and one can view the conversion of the right-hand side of
formula (8.8.43) of \cite{FLM} into the illuminating
generating-function form (\ref{newjacobiiden}) as a nontrivial example
of this philosophy.  (A much more elementary use of generating
functions, central to the present work, is formula (\ref{infinitdil}),
relating a generating function ranging over the number of derivatives
of a formal series with a formal change of variables, and expressing
the fact that $D_x$ is a formal infinitesimal dilation.)

After this review and elaboration of the results in \cite{L2}, below
we present our results: the extension of (\ref{newjacobiiden}) to
expressions of the type $X(Y[u,y]v,x)$, the use of the resulting
identity to compute commutators of such expressions, and the
application of these general principles to the recovery of Theorem 3.1
of \cite{L2} (see Theorem \ref{theoremforLbar} above)---a commutator
formula for certain vertex operators that in turn yields a central
extension of a Lie algebra of differential operators and the phenomena
found by Bloch in \cite{Bl}.

\renewcommand{\theequation}{\thesection.\arabic{equation}}
\renewcommand{\therema}{\thesection.\arabic{rema}}
\setcounter{equation}{0}
\setcounter{rema}{0}

\section{Main results}

We want to bracket expressions of the type $X(Y[u,y]v,x)$ and thereby
to recover Theorem \ref{theoremforLbar}, for the reasons that we have
been discussing.  First, in order to extract the commutator $[X(u,
x_{1}),X(v, x_{2})]$ {}from (\ref{newjacobiiden}), we perform the
usual procedure (recall \cite{FLM}, \cite{FHL}) of extracting the
coefficient of $x_0^{-1}$ (the formal residue ${\rm Res}_{x_0}$ with
respect to the variable $x_0$) on both sides, leaving us with the
desired commutator on the left-hand side.  To compute this formal
residue on the right-hand side, we use the following variant of a
standard formal change-of-variables formula: Let $A$ be a commutative
associative algebra and let $F(x)$ be a formal power series in
$A[[x]]$ with no constant term and such that the coefficient of $x^1$
is invertible in $A$.  Then for $h(x) \in A((x))$,
\begin{equation}\label{res}
{\rm Res}_{x}h(x) = {\rm Res}_{y} (h(F(y))F'(y)).
\end{equation}
This formula holds more generally for $h(x) \in V((x))$, where $V$ is
any $A$-module.  Using this for $A$ the commutative associative
algebra $\Bbb C [x_1, x_1^{-1}]$, $V$ the $A$-module of
operator-valued formal Laurent series in $x_1$ and $x_2$, $F(y) = x_1
(1 - e^y)$, and $h(x_0)$ equal to the left-hand side of
(\ref{newjacobiiden}), an operator-valued formal Laurent series in
$x_0, x_1$ and $x_2$ with only finitely many negative powers of $x_0$,
we obtain:

\begin{theo}\label{commutatorthm}
In the setting of Theorem \ref{newjacobithm},
\begin{equation}\label{newcomm}
[X(u, x_{1}),X(v,x_{2})] = {\rm Res}_{y} \delta
\left(e^{-y}{\displaystyle\frac{x_{1}}{x_{2}}}\right)X(Y[u, y]v,
x_{2}).
\end{equation}
\end{theo}

Here the dummy variable $y$ is actually the variable denoted $y_{01}$
in Theorem \ref{newjacobithm}.

We can generalize Theorem \ref{newjacobithm} by taking the vectors $u$
and $v$ in that result to be appropriate composite expressions, and by
suitably using (\ref{infinitdil}), to obtain:

\begin{theo}\label{generaljacobithm}
In any vertex operator algebra $V$, for any $u_1, v_1, u_2, v_2 \in V$
we have:
\begin{eqnarray}\label{generaljacobiiden}
&x_{0}^{-1}\delta
\left(e^{w_{21}}{\displaystyle\frac{e^{w_1}x_{1}}{x_{0}}}\right)
X(Y[u_{1},y_{1}]v_{1},e^{w_{1}}x_{1})
X(Y[u_{2},y_{2}]v_{2},e^{w_{2}}x_{2}) &\nonumber \\ 
&-x_{0}^{-1}\delta
\left(-e^{w_{12}}{\displaystyle\frac{e^{w_2}x_{2}}{x_{0}}}\right)
X(Y[u_{2},y_{2}]v_{2},e^{w_{2}}x_{2}) 
X(Y[u_{1},y_{1}]v_{1},e^{w_{1}}x_{1})&\nonumber \\ 
&=(e^{w_{2}}x_{2})^{-1} \delta
\left(e^{-w_{01}}{\displaystyle\frac{e^{w_{1}}x_{1}}{e^{w_{2}}x_{2}}}
\right)X(Y[Y[u_{1},y_{1}]v_{1},w_{01}]Y[u_{2},y_{2}]v_{2},e^{w_{2}}x_{2}),&
\end{eqnarray}
where
\begin{equation}\label{wij}
w_{21} = \log \left(1-{\displaystyle\frac{e^{w_{2}}x_{2}}
{e^{w_{1}}x_{1}}}\right),\;\;
w_{12} = \log \left(1-{\displaystyle\frac{e^{w_{1}}x_{1}}
{e^{w_{2}}x_{2}}}\right),\;\;
w_{01} = -\log
\left(1-{\displaystyle\frac{x_{0}}{e^{w_{1}}x_{1}}}\right).
\nonumber
\end{equation}
\end{theo}

As in Theorem \ref{newjacobithm}, we choose to write the
delta-function expressions in this way to exhibit the variables
appropriate to the context.

{}From either Theorem \ref{commutatorthm} or Theorem
\ref{generaljacobithm} we obtain the corresponding commutator formula:

\begin{theo}\label{generalcommutatorthm}
In the setting of Theorem \ref{generaljacobithm},
\begin{eqnarray}\label{generalcommutatoriden}
&[X(Y[u_{1},y_{1}]v_{1},e^{w_{1}}x_{1}),
X(Y[u_{2},y_{2}]v_{2},e^{w_{2}}x_{2})] &\nonumber \\ 
&={\rm Res}_{w} 
\delta\left(e^{-w}{\displaystyle\frac{e^{w_{1}}x_{1}}{e^{w_{2}}x_{2}}}
\right)X(Y[Y[u_{1},y_{1}]v_{1},w]Y[u_{2},y_{2}]v_{2},e^{w_{2}}x_{2}).&
\end{eqnarray}
\end{theo}

The dummy variable $w$ is the same as the variable denoted $w_{01}$ in
Theorem \ref{generaljacobithm}.  Of course, both Theorem
\ref{generaljacobithm} and Theorem \ref{generalcommutatorthm} extend
to multiple expressions of the type in the right-hand side of
(\ref{several}).

These and related general results can be applied to a wide variety of
special situations.  Here we indicate the main steps in the conceptual
recovery of Bloch's formulas, as reformulated and generalized in
Theorem 3.1 of \cite{L2} (recalled in Theorem \ref{theoremforLbar}
above).

We of course work in the special setting discussed in
(\ref{Yh(-1)})--(\ref{hhsimnordplushh}) above.  We can use Theorem
\ref{generalcommutatorthm} to compute
\begin{equation}\label{specialcomm}
[X(Y[v_{0},y_{1}]v_{0},e^{w_{1}}x_{1}),X(Y[v_{0},y_{2}]v_{0},e^{w_{2}}x_{2})]
\end{equation}
(the case $u_1 = v_1 = u_2 = v_2 = v_0$ in
(\ref{generalcommutatoriden})).  It is convenient to take $w_1 = w_2 =
0$ at first, and then to apply formula (\ref{infinitdil}) to restore
$e^{w_1}$ and $e^{w_2}$ at the end.  In order to evaluate and simplify
the right-hand side of (\ref{generalcommutatoriden}), we
systematically use Zhu's theorem (mentioned above) that the Jacobi
identity holds for the operators $Y[u,x]$.  The analysis of the
right-hand side of (\ref{generalcommutatoriden}) is straightforward
and natural; the details are given in \cite{L3}.  Here we state a
result, in the full generality of Theorem \ref{generalcommutatorthm},
that serves to partially evaluate the expression
(\ref{generalcommutatoriden}); without loss of generality, we state it
for the case $w_1 = w_2 = 0$:

\begin{theo}
In the indicated setting, we have
\begin{eqnarray*}
&[X(Y[u_{1},y_{1}]v_{1},x_{1}),X(Y[u_{2},y_{2}]v_{2},x_{2})] &\nonumber \\ 
&={\rm Res}_{t_{1}}e^{t_{1}\frac{\partial}{\partial y_{1}}}
\left(\delta\left(e^{-t_{1}}{\displaystyle\frac{x_{1}}{x_{2}}}\right)
X(Y[u_{2},y_{2}]Y[u_{1},y_{1}]Y[v_{1},t_{1}]v_{2},x_{2})\right)
&\nonumber \\ 
&+{\rm Res}_{t_{2}}e^{-t_{2}\frac{\partial}{\partial y_{1}}}
\left(\delta\left(e^{y_{1}}{\displaystyle\frac{x_{1}}{x_{2}}}\right)
X(Y[u_{2},y_{2}]Y[v_{1},-y_{1}]Y[u_{1},t_{2}]v_{2},x_{2})\right)
&\nonumber \\ 
&-{\rm Res}_{t_{3}}e^{-t_{3}\frac{\partial}{\partial y_{2}}}
\left(\delta\left(e^{-y_{2}}{\displaystyle\frac{x_{1}}{x_{2}}}\right)
X(Y[Y[u_{1},y_{1}]Y[u_{2},t_{3}]v_{1},y_{2}]v_{2},x_{2})\right)
&\nonumber \\ 
&-{\rm Res}_{t_{4}}e^{-t_{4}\frac{\partial}{\partial y_{2}}}
\left(\delta\left(e^{-y_{2}+y_{1}}{\displaystyle\frac{x_{1}}{x_{2}}}\right)
X(Y[Y[Y[u_{2},t_{4}]u_{1},y_{1}]v_{1},y_{2}-y_{1}]v_{2},x_{2})\right).
&
\end{eqnarray*}
\end{theo}

Even in this full generality, the features of Theorem
\ref{theoremforLbar} are already starting to appear here.  We can
apply this result to the special case (\ref{specialcomm}).  As in the
``direct'' proof of Theorem 3.1 of \cite{L2} (using the ingredients of
Section 3 of \cite{L2} rather than the present general
considerations), rather subtle cancellation occurs on the right-hand
side, but here the cancellation occurs in a more conceptual context.
The result is:

\begin{theo}
Applied to the case (\ref{specialcomm}), Theorem
\ref{generalcommutatorthm} yields precisely the assertion of Theorem
3.1 of \cite{L2}.
\end{theo}

We have thus ``explained'' both the zeta-function-modified operators
and their commutators as studied in \cite{Bl}.

{\small \sc Department of Mathematics, Rutgers University, Piscataway,
NJ 08854}

{\em E-mail address}: lepowsky@math.rutgers.edu

\end{document}